\documentstyle[12pt]{article}
\begin{document}

\newcommand\comp{\circ}
\newtheorem{theorem}{Theorem}[section]
\newtheorem{example}[theorem]{Example}
\newtheorem{conjecture}[theorem]{Conjecture}
\newtheorem{examples}[theorem]{Examples}
\newtheorem{proposition}{Proposition}[section]
\newtheorem{fact}[theorem]{Fact}
\newtheorem{problem}[theorem]{Open problem}

\newcommand{\I}{\mbox{{\bf I}}}
\newcommand{\A}{\mbox{{\bf Ax \ }}}
\newcommand{\qed}{\hspace*{\fill}\rule{1 ex}{1.5 ex}\\}

\newtheorem{corollary}[theorem]{Corollary}
\newtheorem{definition}[theorem]{Definition}
\newtheorem{remark}[theorem]{Remark}
\newtheorem{lemma}[theorem]{Lemma}
\newtheorem{claim}{Claim}[theorem]
\newtheorem{construction}[theorem]{Construction}
\newenvironment{proof}{\begin{trivlist}\item[]{\bf
Proof.}}{\qed\end{trivlist}}

\newcommand{\comment}[1]{\typeout{swallowing comment}}

\hyphenation{equi-v-al-ent dimen-s-ion mon-ad-ic}

\title{ On Weak and Strong Interpolation in Algebraic Logics }
\author{G\'abor S\'agi\thanks{Supported by Hungarian
National Foundation for Scientific Research grant D042177.},
Saharon Shelah\footnote{The second author would like to thank the
Israel Science Foundation for partial support of this research
(Grant no. 242/03). Publication 864.}}


\maketitle
\begin{abstract}
We show that there is a restriction, or modification of the
finite-variable fragments of First Order Logic in which a weak
form of Craig's Interpolation Theorem holds but a strong form of
this theorem does not hold. Translating these results into
Algebraic Logic we obtain a finitely axiomatizable subvariety of
finite dimensional Representable Cylindric Algebras that has the
Strong Amalgamation Property but does not have the
Superamalgamation Property. This settles a conjecture of Pigozzi
\cite{pigo}.
\\
\\
{\em AMS Classification: } 03C40, 03G15.  \\
{\em Keywords:} Craig Interpolation, Strong Amalgamation,
Superamalgamation, Varieties of Cylindric Algebras.
\end{abstract}

\section{Introduction}
\label{intro}

\ \indent Formula interpolation in different logics is a classical
and rapidly growing research area. In this note we give a
modification of finite variable fragments of First Order Logic in
which a weak version of Craig's Interpolation Theorem holds but a
strong version of this theorem does not hold. To do this we will
use classical methods and results
of model theory of First Order Logic. \\
\indent A traditional approach for investigating interpolation
properties of logics is to ''algebraize'' the question, that is,
after reformulating semantics in an algebraic way, interpolation,
definability and related problems can also be considered as
properties of the (variety of) algebras obtained by the above
reformulation. Algebras obtained by algebraizing semantics are
called ''meaning algebras''. As it is well known, the meaning
algebras of first order logics are different classes of
representable cylindric algebras ($RCA_{n}$ will denote the class
of $n$ dimensional representable cylindric algebras; in this paper
$n$ will always be a finite number). For more details we refer to
\cite{HMT1} and \cite{HMT2}.
\\
\indent It turned out that interpolation properties on the
logical side correspond amalgamation properties on the algebraic
side (see, for example, \cite{pigo} or Theorem 6.15 of \cite{ans}
and references therein ). Similarly, Beth Definability Property
(on the logical side) corresponds to surjectiveness of the
epimorphisms in the category of meaning alegabras, see \cite{nemeti} or Theorem 6.11 of \cite{ans}. \\
\indent Our results can also be translated into Algebraic Logic.
As we will show in Theorem \ref{mainthm}, for $n \geq 3$, there
is a subvariety $U_{n}$ of
$RCA_{n}$ that has the Strong Amalgamation Property ($SAP$ for
short) but does not have the Superamalgamation Property ($SUPAP$
for short). This settles a conjecture of Pigozzi in \cite{pigo}
(see page 313, Remark 2.1.21 therein). This was a long standing
open problem in Algebraic Logic.
So Theorem \ref{mainthm} can be considered as the main result of
this note. We should mention the following earlier related
results. Comer \cite{comer} proved that if $n \geq 2$ then
$RCA_{n}$ does not have the Amalgamation Property and Maksimova
\cite{maksi} has shown the existence of a BAO-type\footnote{BAO
stands for Boolean Algebra with Operators} variety that has $SAP$
but doesn't have $SUPAP$. The essential difference between this
result and our Theorem \ref{mainthm} is that in our case the
variety $U_{n}$ is a subvariety of $RCA_{n}$, as originally
Pigozzi's question required. Indeed, as shown in Sections
\ref{interpol} and \ref{defin} this has immediate consequences of
interpolation and definability properties for
some modifications of First Order Logic restricted to finitely many variables. \\
\indent In Section \ref{logdef} we give some basic definitions
about logics in general and then introduce and investigate the
model theory of a certain modification of finite variable
fragment of First Order Logic, this fragment will be called ${\cal
U}_{n}$. In Section \ref{interpol} we show that ${\cal U}_{n}$
satisfies a weak version of Craig's Interpolation Theorem, but
doesn't satisfy the strong version of it. In Section \ref{defin}
we show that ${\cal U}_{n}$ satisfies Beth's Theorem on implicit
and explicit definability. Finally, in Section \ref{cyla} we
translate these results into algebraic form and in Theorem
\ref{mainthm} we settle Pigozzi's conjecture: there is a finite
dimensional, finitely axiomatizable subvariety of representable
cylindric algebras, that has $SAP$ but does not have $SUPAP$. \\
\indent We conclude this section by summing up our system of
notation. \\
\indent Every ordinal is the set of smaller ordinals and natural
numbers are identified with finite ordinals. Throughout, $\omega$
denotes the smallest infinite ordinal. If $A$ and $B$ are sets,
then ${}^{A}B$ denotes the set of functions whose domain is $A$
and whose range is a subset of $B$.

\section{A Portion of First Order Logic}
\label{logdef}


\ \indent In this note by a logic we mean a triplet ${\cal J} =
\langle F, K, \models \rangle$ where $F$ is the set of formulas of
${\cal J}$, $K$ is the class of models of ${\cal J}$ and
$\models$ is the satisfaction relation. Often, formulas have a
structure: there are a {\em vocabulary} and a set of rules with
which one can build formulas from elements of the vocabulary.
Strictly speaking, in this case we obtain different logics for
different vocabularies. Sometimes these families of logics have
been regarded as a pair ${\cal J} = \langle S, K \rangle$ where
$S$ is a function on vocabularies associating the set of formulas
$F=S(V)$ of ${\cal J}$ and the satisfaction relation
$\models_{S(V)}$ of ${\cal J}$ to a given vocabulary $V$. \\
\indent Particularly, when one deals with a concrete first order
language, one should specify the names and arities of relation
and function symbols to be used. Such a specification will also
be called a vocabulary (for first order languages). Throughout
this paper we will deal with variants, modifications and portions
of First Order Logic. For a given vocabulary $V$, $F=F_{V}$ will
be the set of formulas of First Order Logic restricted to
individual variables $\{ v_{0},...,v_{n-1} \}$ ($n$ is fixed and
finite). In addition, in this note, the class $K$ of models of a
logic will always be a subclass of ordinary relational structures
and the
satisfaction relation will be the same as in ordinary First Order Logic. \\
\indent If $\varphi$ is a formula of such a logic ${\cal J}$ then
$voc(\varphi)$ denotes the smallest vocabulary for which
$\varphi$ is really a first order formula. If $V$ is a vocabulary
then ${\cal J}[V]$ denotes the logic in which the set of formulas
consists of formulas of ${\cal J}$ whose vocabularies are
contained in $V$, the class of models of ${\cal J}[V]$ is the
class of $V$-reducts of models of ${\cal J}$ and the satisfaction
relation of ${\cal J}[V]$ is the same as that of ${\cal J}$.
Similarly, if $V \subseteq W$ are vocabularies and ${\cal A}$ is
a model for ${\cal J}$ with vocabulary $W$ then ${\cal A}|_{V}$
denotes the reduct of ${\cal A}$ in which only elements of $V$
interpreted as basic relations (or functions). In
this note we will deal with vocabularies consisting relation symbols only. \\
\indent Truth, meaning, and semantical consequence defined the
obvious way, that is, these notions simply inherited from the
first order case. Similarly, some concepts, methods, etc. of First
Order Logic (like isomorphism, elementary equivalence,
(generated) submodels of a structure) will be used the obvious way
without any additional explanation. \\
\indent As we mentioned, throughout the paper $n \in \omega$ is a
fixed natural number. ${\cal L}_{n}$ denotes usual First Order
Logic restricted to the first $n$ individual variables.

\begin{definition}
Let $A$ be a non-empty set, let $k \in \omega$ and let $\bar{s} \in {}^{k}A$. Then \\
\\
\centerline{$ker(\bar{s}) = \{ \langle i,j \rangle \in {}^{2}k:
s_{i}=s_{j} \}$.} \\
\\
If $U_{0}
\subseteq A$ and $\bar{z} \in {}^{k}A$ then $\bar{s} \sim_{A,U_{0},k} \bar{z}$ means that \\
\indent (i) $ker(\bar{s}) = ker(\bar{z})$ and \\
\indent (ii) $(\forall i \in k)[ s_{i} \in U_{0} \Leftrightarrow
z_{i} \in U_{0}]$. \\
Sometimes we will simply write $\sim_{k}$ or $\sim$ in place of
$\sim_{A,U_{0},k}$.
\end{definition}

\begin{definition}
\label{ustrdef} A relational structure ${\cal A} = \langle A,
U_{0}, R_{i} \rangle_{i \in V}$ is defined to be an $U$-structure
(for the vocabulary $V=V[{\cal A}]$), if \\
\indent $\bullet$ $U_{0} \subseteq A$, \\
\indent $\bullet$ $|U_{0}| \geq n, |A-U_{0}| \geq n$ and \\
\indent $\bullet$ for any $i \in V$, if $R_{i}$ is $k$-ary,
$\bar{s} \in R_{i}, \bar{s} \sim_{k} \bar{z}$
then $\bar{z} \in R_{i}$. \\
$A$ is the universe of ${\cal A}$ and $U_{0}$ will be called the
core of ${\cal A}$.
\end{definition}

In Sections \ref{logdef}, \ref{interpol} and \ref{defin}
$U$-structures have been treated as special first order
relational structures, that is, every relation has a finite
arity, these arities may be different for different relations. \\
\indent Let ${\cal A}$ be an $U$-structure with core $U_{0}$. It
is easy to see that a permutation of the universe of ${\cal A}$
mapping $U_{0}$ onto itself is an automorphism of ${\cal
A}$. \\
\indent It should be emphasized, that the core of an
$U$-structure ${\cal A}$ is not a basic relation of ${\cal A}$,
that is, the core relation a priori doesn't have a name in the
vocabulary of ${\cal A}$. Sometimes the core may be defined
somehow, in some other $U$-structures the core cannot be defined
by first order formulas. The core relation provides some extra
structure for $U$-structures which will be used extensively
below. The core of an $U$-structure ${\cal A}$ will be denoted by
$U_{0}^{\cal A}$ or simply by $U_{0}$ when ${\cal A}$ is clear
from the context.
\\
\indent Throughout this paper by a ''definable relation'' we
mean a relation which is definable by a formula of ${\cal L}_{n}[V]$ without parameters
(if the vocabulary $V$ is clear from the context, we omit it). \\
\indent If $A$ is any set then $A^{*n} = \{ s \in {}^{n}A:
(\forall i\not=j \in n) s_{i} \not= s_{j} \}$. Clearly, this
relation is definable from the identity relation.
In order to keep notation simpler, we will identify this relation
by one of it's defining formulas and sometimes we will write
''$A^{*n}$'' in the middle of another
formula. \\

\begin{definition}
\label{strongustrdef} An $U$-structure ${\cal A}$ is defined to be
a strong $U$-structure if the following holds. If $V$ is any
sub-vocabulary
of the vocabulary of ${\cal A}$ such that \\
\indent $\bullet$ $U_{0}^{\cal A}$ is not ${\cal L}_{n}$-definable
in ${\cal A}|_{V}$ and \\
\indent $\bullet$ $R$ is an ${\cal L}_{n}$-definable $m$-ary
relation of ${\cal A}|_{V}$ (for some $m \leq n$) such that $R
\subseteq A^{*m}$ and $i,j \in m, \ i \not= j$ then \\
${\cal A} \models R \Leftrightarrow (\exists v_{i}R \wedge
\exists v_{j} R) \wedge A^{*m}$ (more precisely, letting $\bar{v}
= \langle v_{0},...,v_{n-1} \rangle$ we require ${\cal A} \models
(\forall \bar{v}) [R(\bar{v}) \Leftrightarrow \exists v_{i}
R(\bar{v}) \wedge \exists v_{j}R(\bar{v}) \wedge
A^{*m}(\bar{v})]$).
\end{definition}

${\cal U}_{n}$ will denote the logic in which the set of formulas
is the same as in ${\cal L}_{n}$ and the class of models of
${\cal U}_{n}$ is the class of strong $U$-structures. We will say
that a relation (in an arbitrary structure) is ${\cal
U}_{n}$-definable iff it is definable by a formula of ${\cal
U}_{n}$. Similarly, two relational structures are called ${\cal
U}_{n}$-elementarily equivalent iff they satisfy the same
formulas of ${\cal U}_{n}$. \\
\\
\indent In the previous definitions ''U'' stand for
''unary-generated'', this choice of naming will be explained in
Section \ref{cyla} below. According to the previous definition,
the notion of strong $U$-structures depends on $n$, therefore
strictly speaking, instead of ''strong $U$-structure'' we should
write ''strong $U$-structure for some $n$''. For simplicity we
don't indicate
$n$; it will always be clear from the context. \\
\indent We call the attention that in ${\cal U}_{n}$ function
symbols are not part of the vocabulary, that is, all the
vocabularies contain relation symbols only. \\
\indent It is easy to see that strong $U$-structures exist. We
will show this in Theorem \ref{uexistthm} below.

\begin{lemma}
\label{basiclemma} Suppose ${\cal A}$ is an $U$-structure with
universe $A$ and core $U_{0}$. \\
\indent (1) If $R$ is a definable unary relation of ${\cal A}$
then $R \in \{ \emptyset, A, U_{0}, A-U_{0} \}$. Thus, at
most four definable unary relations exist in an $U$-structure. \\
\indent (2) Let ${\cal A}' = \langle A, U_{0} \rangle$ be the
structure whose universe is the same as that of ${\cal A}$ and
whose unique basic relation is $U_{0}$. If $R$ is a ${\cal
U}_{n}$-definable relation of
${\cal A}$ then $R$ is definable in ${\cal A}'$ as well. \\
\indent (3) If $R$ is a definable relation in ${\cal A}$, $\bar{s}
\in R$ and $\bar{z} \sim \bar{s}$ then $\bar{z} \in R$.
\end{lemma}

\begin{proof}
First observe the following. If $f$ is any permutation of $A$
preserving $U_{0}$ (i.e. mapping it onto itself) then for any $k
\in \omega$ and $\bar{s} \in {}^{k}A$ we have $\bar{s} \sim_{k}
f(\bar{s})$. Therefore by Definition \ref{ustrdef} $f$ is an
automorphism of ${\cal A}$. Now suppose $a \in U_{0} \cap R$ and
$b \in U_{0}$. Then there is an automorphism $f$ of ${\cal A}$
mapping $a$ onto $b$. Since $R$ is definable, $f$ preserves $R$,
thus $b \in R$. It follows, that if $R \cap U_{0} \not=
\emptyset$ then $U_{0} \subseteq R$. Similarly, if $R \cap
(A-U_{0}) \not= \emptyset$ then $A-U_{0}
\subseteq R$, whence (1) follows. \\
\indent Now we turn to prove (2). Suppose that the arity of $R$
is $k$. Since in ${\cal U}_{n}$ there
are only $n$ individual variables, it follows that $k \leq n$. \\
\indent For any equivalence relation $e \subseteq {}^{2}k$ let
$D_{e} = \{s \in {}^{k}A: ker(s) = e \}$. For any $f: k/e
\rightarrow 2$ let $At(f) = \{s \in D_{e}: (\forall i \in k) s_{i}
\in U_{0} \Leftrightarrow f(i/e) = 0 \}$. Clearly, every $D_{e}$
is definable in ${\cal A}'$ (in fact, these relations are
definable from the identity (equality) relation with a
quantifier-free formula of ${\cal U}_{n}$ which doesn't contain
any other basic relation symbol). Similarly, every $At(f)$ is
definable in ${\cal A}'$. Since ${}^{k}A$ is the disjoint union of
the $D_{e}$'s, it is enough to show that for all $e$ the relation
$R \cap D_{e}$ is definable in ${\cal A}'$. Let $e$ be fixed. If
$R \cap D_{e} = \emptyset$ then it is definable in
${\cal A}'$, so we may assume $R \cap D_{e} \not= \emptyset$. \\
\indent Suppose $\bar{s} \in R \cap D_{e} \cap At(f)$ for some $f:
k/e \rightarrow 2$. We claim that in this case $At(f) \subseteq R
\cap D_{e}$. To check this suppose $\bar{z} \in At(f)$. Let $g$ be
the partial function on $A$ mapping each $s_{i}$ onto $z_{i}$.
Since $ker(\bar{z}) = e = ker(\bar{s})$, $g$ is a well defined
partial function and moreover $g$ is injective. In addition, for
every $i \in k$, $s_{i} \in U_{0} \Leftrightarrow z_{i} \in
U_{0}$. Therefore there is a permutation $h$ of $A$ extending $g$
and preserving $U_{0}$. As observed at the beginning of the proof
of (1), $h$ is an automorphism of ${\cal A}$. Since $R \cap D_{e}$
is definable in ${\cal A}$, it follows that $h$ preserves
$R \cap D_{e}$ and thus $\bar{z} = f(\bar{s}) \in R \cap D_{e}$, as desired. \\
\indent Now let $S = \{ At(f): At(f) \cap R \cap D_{e}
\not=\emptyset,$ \
$f \in {}^{k/e}2 \}$ and let $P = \cup S$. Clearly, $P$ is
definable in ${\cal A}'$. We claim, that $P = R \cap D_{e}$. By
the previous paragraph we have $P \subseteq R \cap D_{e}$. On the
other hand, if $\bar{s} \in R \cap D_{e}$, then for the function
$f: k / e \rightarrow 2$, $f(i/e) = 0 \Leftrightarrow s_{i} \in
U_{0}$ we have $\bar{s} \in At(f)$, therefore every element of $R
\cap D_{e}$ is contained in an element of $S$
and thus $R \cap D_{e} \subseteq P$. \\
\indent For (3) observe that for any $At(f)$, if $\bar{s} \sim
\bar{z}$ and $\bar{s} \in At(f)$ then $\bar{z} \in At(f)$. Now by
the previous proof of (2), if $\bar{s} \in R$ then $\bar{s} \in
At(f) \subseteq R$ for some $f$ and therefore $\bar{z} \in At(f)$
whence $\bar{z} \in R$.
\end{proof}


Below we will associate {\em Cylindric Set Algebras} with
relational structures in the usual way. For
completeness we recall here the details. \\
\indent suppose ${\cal A}$ is a relational structure. It's
$n$-dimensional Cylindric Set Algebra will be denoted by
$Cs_{n}({\cal A})$. Roughly speaking, the elements of
$Cs_{n}({\cal A})$ are the ${\cal L}_{n}$-definable relations of
${\cal A}$. To be more precise, elements of $Cs_{n}({\cal A})$
are $n$-ary relations. If $\varphi(v_{0},...,v_{m-1})$ is a
formula of ${\cal L}_{n}$ in the vocabulary of ${\cal A}$ with
free variables as indicated then $\varphi$ defines an $m$-ary
relation in ${\cal A}$. The corresponding element of
$Cs_{n}({\cal A})$ is the $n$-ary relation $[\varphi] = \{
\bar{s} \in {}^{n}A: {\cal A} \models \varphi[\bar{s}] \}$.
The $n$-dimensional Cylindric Set Algebra $Cs_{n}({\cal A})$ of
${\cal A}$ is the following algebra ${\cal B} = \langle X;
\cap,-,C_{i},D_{i,j} \rangle_{i,j \in n}$. Here $X=\{ [\varphi]:
\varphi$ is a formula of ${\cal L}_{n} \}$ is the set of elements
of ${\cal B}$. The operations $\cap$ and $-$ are set-theoretic
intersection and complementation (w.r.t. ${}^{n}A$),
respectively. Then for any formulas $\varphi, \psi$ of ${\cal L}_{n}$ one has \\
\\
\indent $[\varphi] \cap [\psi] = [\varphi \wedge \psi]$ and \\
\indent $-[\varphi] = [\neg \varphi]$.\\
\\
In addition, $C_{i}$ is a unary operation and $D_{i,j}$ is a
0-ary operation for every $i,j \in n$. These operations
correspond to the semantics of existential quantifier and the
equality symbol of First Order Logic. In more detail, if
$[\varphi] \in X$ is an element of
${\cal B}$
then \\
\\
\indent $C_{i}([\varphi]) = [\exists v_{i} \varphi]$ and \\
\indent $D_{i,j} = \{ s \in {}^{n}A: s_{i} = s_{j} \}$. \\
\\
Let $V$ be the vocabulary of ${\cal A}$. Clearly, $\{[R_{i}]: i
\in V \}$ is a set of generators of $Cs_{n}({\cal A})$.

\begin{lemma}
\label{upclosed} The class of (strong) $U$-structures is closed
under ultraproducts.
\end{lemma}

\begin{proof}
Let $\langle {\cal A}_{i}: i \in I \rangle$ be a system of
$U$-structures and let ${\cal F}$ be an ultrafilter on $I$. Then
$\Pi_{i \in I} {\cal A}_{i} / {\cal F}$ is an $U$-structure (with
core $\Pi_{i \in I} U_{0}^{{\cal A}_{i}} / {\cal F}$) because the
requirements of Definition \ref{ustrdef} can be expressed by
first order formulas in the expanded vocabulary in
which there is an extra symbol for the core relation. \\
\indent Now let $\langle {\cal A}_{i}: i \in I \rangle$ be a
system of strong $U$-structures, let ${\cal F}$ be an ultrafilter
on $I$ and let ${\cal A} = \Pi_{i \in I}{\cal A}_{i} / {\cal F}$.
Let $V=\{ R_{0},...,R_{h} \}$ be a finite
sub-vocabulary of the common vocabulary of the previous system of structures. \\
\indent Let $J = \{ i \in I: U_{0}^{{\cal A}_{i}}$ is definable
in ${\cal A}_{i}|_{V} \}$. We will show that $J \in {\cal F}$
implies that the core $\Pi_{i \in I}U_{0}^{{\cal A}_{i}} / {\cal
F}$ of ${\cal A}$ is definable in ${\cal A}|_{V}$. So suppose $J
\in {\cal F}$. For each $i \in J$ fix a ${\cal U}_{n}[V]$-formula
$\varphi_{i}$ defining $U_{0}^{{\cal A}_{i}}$ in ${\cal
A}_{i}|_{V}$.
Observe, that there is a finite number $N_{0}$ such that for all
$k \leq n$, for all $A$ and for all $U_{0}$ the equivalence
relation $\sim_{A,U_{0},k}$ has at most $N_{0}$ equivalence
classes. Hence, by Lemma \ref{basiclemma} (3) there is a finite
number $N_{1}$ such that $|Cs_{n}({\cal A}_{i}|_{V})| \leq N_{1}$
for all $i \in I$. For any $i \in I$ let $Cs_{n}^{+}({\cal
A}_{i}|_{V}) = \langle Cs_{n}({\cal A}_{i}|_{V}), [R_{i}]
\rangle_{i \in V}$, that is, $Cs_{n}^{+}({\cal A}_{i}|_{V})$ is
$Cs_{n}({\cal A}_{i}|_{V})$ expanded with the relations
corresponding to the interpretations of elements of $V$. The set
$\{ [R_{i}]: i \in V \}$ generates $Cs_{n}({\cal A}_{i}|_{V})$
therefore for each $i \in I$ there is a finite set $T_{i}$ of
cylindric terms such that for every $a \in Cs_{n}({\cal
A}_{i}|_{V})$ there is a $t \in T_{i}$ with
$a=t([R_{0}],...,[R_{h}])$. We may assume that $Cs_{n}^{+}({\cal
A}_{i}|_{V}) \cong Cs_{n}^{+}({\cal A}_{j}|_{V})$ implies $T_{i}
= T_{j}$ for all $i,j \in I$. Since $V$ is finite and
$|Cs_{n}({\cal A}_{i}|_{V})| \leq N_{1}$ for all $i \in I$ there
exist a finite set $T$ of cylindric terms and $K \subseteq J$
such that $K \in {\cal F}$ and for every $i,j \in K$ we have
$Cs_{n}^{+}({\cal A}_{i}|_{V}) \cong Cs_{n}^{+}({\cal
A}_{j}|_{V})$ and $T = T_{i}$. Hence there are a $t \in T$ and $L
\subseteq K$ such that $L \in {\cal F}$ and for every $i \in L$
we have $[\varphi_{i}] = t([R_{0}],...,[R_{h}])$. Let $\varphi$ be
the formula corresponding to $t([R_{0}],...,[R_{h}])$. Then
clearly,
$\varphi$ defines the core of ${\cal A}$. \\
\indent Next we show that ${\cal A}$ is a strong $U$-structure.
Suppose $W$ is (an arbitrary, not necessarily finite)
sub-vocabulary of the vocabulary of ${\cal A}$ such that the core
of ${\cal A}$ is not definable in ${\cal A}|_{W}$ and $R$ is an
${\cal L}_{n}$-definable relation of ${\cal A}|_{W}$.  Then there
is a finite sub-vocabulary $V \subseteq W$ such that $R$ is
definable in ${\cal A}|_{V}$ and still, the core of ${\cal A}$ is
not definable in ${\cal A}|_{V}$. Applying the result of the
previous paragraph to this $V$, it follows that $J \not\in {\cal
F}$.
Finally observe that the required property of $R$ (described in
Definition \ref{strongustrdef}) can be expressed by a formula of
${\cal L}_{n}[V]$ and this formula is true in ${\cal A}$ since for
every $i \in I$ ${\cal A}_{i}$ is a strong $U$-structure.
\end{proof}

\begin{theorem}
\label{uexistthm} (1) There exists a strong $U$-structure. \\
\indent (2) There exists a structure ${\cal A} = \langle A,
U_{0},P,Q \rangle$ which is a strong $U$-structure with core
$U_{0}$ such that $P$ and $Q$ are unary relations and $P=Q=U_{0}$.
\end{theorem}

\begin{proof}
Since (2) implies (1), it is enough to prove (2). Let $A$ be any
countably infinite set, let $U_{0} \subseteq A$ be such that
$|U_{0}| = |A-U_{0}| = \aleph_{0}$ and finally let $P=Q=U_{0}$. We
have to show that ${\cal A} = \langle A,U_{0},P,Q \rangle$ is a
strong $U$-structure. Since $U_{0}$ is infinite and the basic
relations of ${\cal A}$ are unary, ${\cal A}$ satisfies
Definition \ref{ustrdef} for every $n \in \omega$. Thus, ${\cal
A}$ is an $U$-structure (for any $n \in \omega$). \\
\indent Now suppose $V$ is a sub-vocabulary of the vocabulary of
${\cal A}$ such that $U_{0}$ is not ${\cal U}_{n}$-definable in
${\cal A}|_{V}$. It follows that $V$ contains the equality symbol
only. Suppose $R$ is an $m$-ary relation ${\cal U}_{n}$-definable
in ${\cal A}|_{V}$ such that $R \subseteq A^{*m}$. If $R=
\emptyset$ then Definition \ref{strongustrdef} holds for $R$. Now
suppose $\bar{s} \in R$ and $\bar{z} \in A^{*m}$. Then there is a
permutation $f$ of $A$ mapping $\bar{s}$ onto $\bar{z}$. Since
permutations preserve the identity relation and $R$ is definable
in ${\cal A}|_{V}$, it follows that $f$ preserves $R$ and
therefore $\bar{z} \in R$. Since $\bar{z} \in A^{*m}$ was
arbitrary, $R = A^{*m}$. Clearly, this relation satisfies the
requirements of Definition \ref{strongustrdef}. So ${\cal A}$ is
a strong $U$-structure, as desired.
\end{proof}

\begin{theorem}
\label{autothm} Suppose ${\cal A}$ is a strong $U$-structure. If
$\bar{a}, \bar{b} \in A$ satisfy the same ${\cal U}_{n}$-formulas
in ${\cal A}$ then there is an automorphism of ${\cal A}$ mapping
$\bar{a}$ onto $\bar{b}$.
\end{theorem}

\begin{proof}
Let $U_{0}$ be the core of ${\cal A}$. First suppose that $U_{0}$
can be defined in ${\cal A}$ by a ${\cal U}_{n}$-formula. In this
case (since $\bar{a}$ and $\bar{b}$ satisfy the same ${\cal
U}_{n}$-formulas in ${\cal A}$) we have $\bar{a} \sim \bar{b}$.
Then there is a permutation $f$ of $A$ preserving $U_{0}$ and
mapping $\bar{a}$ onto $\bar{b}$. Then $f$ is an automorphism of
${\cal A}'=\langle A,U_{0} \rangle$ hence it also preserves all
the relations definable in ${\cal A}'$. Hence by Lemma
\ref{basiclemma} (2) $f$ preserves every definable relation of
${\cal A}$ as well, particularly, $f$ is an automorphism of ${\cal
A}$. (There is another way to prove that $f$ is an automorphism of
${\cal A}$: since $f$ preserves $U_{0}$, for every tuple $\bar{s}
\in A$ we have $\bar{s} \sim f(\bar{s})$ hence by Lemma
\ref{basiclemma} (3) it also follows that $f$ is an automorphism
of ${\cal A}$.) \\
\indent Now suppose $U_{0}$ is not ${\cal U}_{n}$-definable in
${\cal A}$. We claim that every relation $R$ definable in ${\cal
A}$ is definable using the identity relation only. This will be
proved by induction on the arity of $R$. If $R$ is unary then by
Lemma \ref{basiclemma} (1) $R$ is either the empty set or $R =
A$; in both cases $R$ is ${\cal U}_{n}$-definable from the
identity relation. Now suppose that $k<n$, $R$ is $k+1$-ary, and
that the claim is true for any relation with arity at most $k$.
Again, for any equivalence relation $e \subseteq {}^{2}(k+1)$ let
$D_{e} = \{s \in {}^{k+1}A: ker(s) = e \}$. Clearly, $D_{e}$ is
${\cal U}_{n}$-definable from the identity relation for any $e$
and $R=\cup_{e}(R \cap D_{e})$. Therefore it is enough to show
that $R \cap D_{e}$ is ${\cal U}_{n}$-definable from the identity
relation. Let $m \subseteq k+1$ be a set of representatives for
$e$ and for any $s \in A^{*m}$ let $s' \in {}^{k+1}A$ be the
sequence for which $ker(s') = e$ and $s=s'|m$. Let $Q = \{s \in
A^{*m}: s' \in R \cap D_{e} \}$. Then $Q$ is ${\cal
U}_{n}$-definable and $R \cap D_{e}$ is definable from $Q$ and
from the identity relation. If $Q$ is at most unary then we are
done because of the basic step of the induction. Otherwise there
are distinct $i,j \in m$ and since ${\cal A}$ is a strong
$U$-structure, we have ${\cal A} \models Q \Leftrightarrow
\exists v_{i} Q \wedge \exists v_{j} Q \wedge A^{*m}$. But the
first two relations in the right hand side are at most $k$-ary,
therefore by the induction hypothesis they are ${\cal
U}_{n}$-definable from the identity relation. Hence $Q$ and
therefore $R \cap D_{e}$ is ${\cal
U}_{n}$-definable in the same way, as well. \\
\indent So suppose $U_{0}$ is not ${\cal U}_{n}$-definable in
${\cal A}$ and $\bar{a}$ and $\bar{b}$ satisfy the same ${\cal
U}_{n}$-formulas in ${\cal A}$. Then $ker(\bar{a}) =
ker(\bar{b})$ and hence there is a permutation $f$ of $A$ mapping
$\bar{a}$ onto $\bar{b}$. Therefore $f$ preserves the identity
relation of ${\cal A}$ and thus, by the previous paragraph, $f$
preserves all the definable relations of ${\cal A}$. So $f$ is
the required automorphism of ${\cal A}$.
\end{proof}

Suppose ${\cal B}$ is a substructure of ${\cal A}$. If
$\bar{k},\bar{k}'$ are tuples of $A$ with the same length such
that $k_{j} = k_{j}'$ for every $j \not=i$ then we will write $k
\stackrel{i}{\cong} k'$. Recall that by the Tarski-Vaught test
${\cal B}$ is an elementary substructure of ${\cal A}$ if for any
first order formula $\varphi$ and tuple $\bar{k} \in B$ we have
${\cal A} \models \exists v_{i} \varphi[\bar{k}]$ if and only if
there is another tuple $\bar{k}' \in B$ such that ${\cal A}
\models \varphi[\bar{k}']$ and $\bar{k} \stackrel{i}{\cong}
\bar{k}'$. It is also easy to check that ${\cal B}$ is a ${\cal
U}_{n}$-elementary substructure of ${\cal A}$ if the previous
condition holds for every ${\cal U}_{n}$-formula $\varphi$.

\begin{theorem}
\label{equivthm} (1) Let ${\cal A}$ be a $U$-structure with core
$U_{0}$ and suppose $V \subseteq A$ is such that $|V \cap U_{0}|,
|V - U_{0}| \geq n$. Let ${\cal B}$ be the substructure of ${\cal
A}$ generated by $V$. Then ${\cal B}$ is an $U$-structure (with
core $U_{0} \cap V$) which is a ${\cal U}_{n}$-elementary
substructure of ${\cal A}$.
${\cal A}$ is a strong $U$-structure if and only if so is ${\cal B}$. \\
\indent (2) Suppose ${\cal A}$ and ${\cal B}$ are ${\cal
U}_{n}$-elementarily equivalent $U$-structures with cores
$U_{0},V_{0}$, respectively.
Then any bijection $f:A \rightarrow B$
mapping $U_{0}$ onto $V_{0}$ is an isomorphism between ${\cal A}$
and ${\cal B}$.
\end{theorem}

\begin{proof}
To prove (1) we have to verify (the above recalled version of) the
Tarski-Vaught test. To do this suppose $\bar{k} \in V$ and
$\varphi$ is a ${\cal U}_{n}$-formula such that ${\cal A} \models
\exists v_{i} \varphi[\bar{k}]$. Let $\bar{k}' \in A$ be a tuple
for which $\bar{k} \stackrel{i}{\cong} \bar{k}'$ and ${\cal A}
\models \varphi[\bar{k}']$. By the condition on $V$, there is
another tuple $\bar{h} \in V$ such that $\bar{h} \sim \bar{k}'$
and $\bar{h} \stackrel{i}{\cong} \bar{k}'$. Therefore it follows
from Lemma \ref{basiclemma} (3) that ${\cal A} \models \varphi
[\bar{h}]$. This shows that ${\cal B}$ is a ${\cal
U}_{n}$-elementary substructure of ${\cal A}$. We claim that
${\cal B}$ is an $U$-structure with core $U_{0} \cap V$. To check
this suppose $R$ is an $m$-ary basic relation of ${\cal B}$,
$\bar{s} \sim \bar{z} \in {}^{m}V$ and $\bar{s} \in R^{\cal B}$.
Then $\bar{s} \in R^{\cal A}$ and since ${\cal A}$ is an
$U$-structure, $\bar{z} \in R^{\cal A}$ hence $\bar{z} \in
R^{\cal B}$, as desired.
\\
\indent Now suppose ${\cal A}$ is a strong $U$-structure. First
observe that if $W$ is a sub-vocabulary of the vocabulary of
${\cal A}$ then by elementarity, if the core of ${\cal A}$ is
${\cal U}_{n}$-definable in ${\cal A}|_{W}$ then the core of
${\cal B}$ is also ${\cal U}_{n}$-definable in ${\cal B}|_{W}$.
In addition, if the core of ${\cal A}|_{W}$ is not definable then
by Lemma \ref{basiclemma} (1) the only unary relations definable
in ${\cal A}|_{W}$ are the empty set and the whole universe of
${\cal A}$; thus, the same is true for ${\cal B}|_{W}$ and
therefore in this case the core of ${\cal B}|_{W}$ is also not
${\cal U}_{n}$-definable. Now suppose $W$ is such a sub-vocabulary
that the core of ${\cal B}$ is not ${\cal U}_{n}$-definable in
${\cal B}|_{W}$ and $R^{\cal B}$ is an $m$-ary ${\cal
U}_{n}$-definable relation in ${\cal B}|_{W}$ such that $R^{\cal
B} \subseteq V^{*m}$ and $i,j \in m, \ i \not= j$. Then by
elementarity $R^{\cal A} \subseteq A^{*m}$ because this property
of $R$ can be described by a ${\cal U}_{n}$-formula. As observed,
the core of ${\cal A}$ cannot be defined in ${\cal A}|_{W}$.
Therefore, since ${\cal A}$ is a strong $U$-structure, ${\cal A}
\models R \Leftrightarrow \exists v_{i} R \wedge \exists v_{j} R
\wedge A^{*m}$. Again by elementarity
the same formula is valid in ${\cal B}$, hence ${\cal B}$ is
indeed a strong $U$-structure. A similar argument shows that if
${\cal A}$ is not a strong $U$-structure then ${\cal B}$
is also not a strong $U$-structure. \\
\indent To show (2) let $f: A \rightarrow B$ be any bijection
mapping $U_{0}$ onto $V_{0}$. Then clearly, $f$ is an isomorphism
between $\langle A, U_{0} \rangle $ and $\langle B, V_{0}
\rangle$. Therefore $f$ preserves any relation which can be
defined by a ${\cal U}_{n}$-formula from $U_{0}$. By Lemma
\ref{basiclemma} (2) every definable (particularly every basic)
relation of ${\cal A}$ can be defined from $U_{0}$, thus $f$
preserves them.
\end{proof}

The following is an adaptation of Corollary 6.1.17 of \cite{chk}.

\begin{theorem}
\label{septhm} (Separation Theorem.) \\
Suppose $K_{0}$ and $K_{1}$ are disjoint classes of strong
$U$-structures with same vocabularies such that both $K_{0}$ and
$K_{1}$ are closed under ultraproducts and ${\cal
U}_{n}$-elementary equivalence. Then there is a ${\cal
U}_{n}$-formula $\varphi$ with $K_{0} \models \varphi$ and $K_{1}
\models \neg \varphi$.
\end{theorem}

\begin{proof}
Recall that by Lemma \ref{upclosed} any ultraproduct of strong
$U$-structures is a strong $U$-structure. \\
\indent Let $\Sigma$ be the set of ${\cal U}_{n}$-formulas valid
in $K_{0}$. Suppose, seeking a contradiction, that there is no
$\varphi$ satisfying the requirements of the theorem. It follows
that every finite subset of $\Sigma$ also has a model in $K_{1}$.
Since $K_{1}$ is closed under ultraproducts, there is a strong
$U$-structure ${\cal A}_{1} \in K_{1}$ such that ${\cal A}_{1}
\models \Sigma$. In addition, if $\Psi$ is a finite set of ${\cal
U}_{n}$-formulas valid in ${\cal A}_{1}$ then $\Psi$ has a model
in $K_{0}$ (otherwise $K_{0} \models \neg (\bigwedge \Psi)$ and
hence $ \neg(\bigwedge \Psi) \in \Sigma$ would follow, therefore
we would have ${\cal A}_{1} \models \neg(\bigwedge \Psi)$). Since
$K_{0}$ is closed under ultraproducts there is an ${\cal A}_{0}
\in K_{0}$ which is ${\cal U}_{n}$-elementarily equivalent with
${\cal
A}_{1}$. \\
\indent Summing up, ${\cal A}_{0} \in K_{0}, {\cal A}_{1} \in
K_{1}$ and ${\cal A}_{0}$ and ${\cal A}_{1}$ are ${\cal
U}_{n}$-elementarily equivalent. This is impossible because
$K_{0}$ and $K_{1}$ are disjoint classes and both are closed
under ${\cal U}_{n}$-elementary equivalence.
\end{proof}

\section{Interpolation}
\label{interpol}

\ \indent We start this section by recalling the weak and strong
forms of Craig's Interpolation Theorem. Suppose ${\cal L}$ is a
logic (in the sense of the beginning of Section \ref{logdef}) and
$\varphi$ is a formula of ${\cal L}$.
%
%
Then $\models \varphi$ means that for any model ${\cal A}$ for
${\cal L}$, $\varphi$ is valid in ${\cal A}$.  If $\psi$ is
another formula of ${\cal L}$ then, as expected, $\varphi \models
\psi$ means that $\psi$ is valid in every model in which
$\varphi$ is valid.

\begin{definition}
A logic ${\cal L}$ has the Strong Craig Interpolation Property if
for any pair of formulas $\varphi, \psi$ of ${\cal L}$ the
following holds. If $\models \varphi \Rightarrow \psi$
then there is a formula $\vartheta$ such that $\models (\varphi
\Rightarrow \vartheta) \wedge (\vartheta \Rightarrow \psi)$ and
the relation symbols
occurring in $\vartheta$ occur both in $\varphi$ and in $\psi$. \\
\indent A logic ${\cal L}$ has the Weak Craig Interpolation
Property if for any pair of formulas $\varphi$ and $\psi$ of
${\cal L}$ the following holds. If $\varphi \models \psi$ then
there is a formula $\vartheta$ such that $\varphi \models
\vartheta$ and $\vartheta \models \psi$ and the relation symbols
occurring in $\vartheta$ occur both in $\varphi$ and in $\psi$.
\end{definition}

\begin{lemma}
\label{elemclosedlemma} Suppose $\varphi$ is a ${\cal
U}_{n}$-formula and $V \subseteq voc(\varphi)$ is a vocabulary.
Then the class $K$ of $V$-reducts of ${\cal U}_{n}$-models of
$\varphi$ is closed under ${\cal U}_{n}$-elementary equivalence.
\end{lemma}

\begin{proof}
Suppose ${\cal A}_{0} \in K$ and ${\cal A}_{1}$ is ${\cal U}_{n}$-
elementarily equivalent with ${\cal A}_{0}$. Let $U_{0}$ and
$U_{1}$ be the cores of ${\cal A}_{0}$ and ${\cal A}_{1}$,
respectively. Let ${\cal A}_{0}^{+}$ be an expansion of ${\cal
A}_{0}$ which is a model of $\varphi$. Let $C$ and $C_{0}
\subseteq C$ be sets such that $|C_{0}| \geq |U_{0}|,|U_{1}|$ and
$|C - C_{0}| \geq |A_{0}-U_{0}|, |A_{1} - U_{1}|$. According to
these cardinal conditions we may (and will) assume $U_{0}, U_{1}
\subseteq
C_{0}$ and $A_{0}-U_{0}, A_{1}-U_{1} \subseteq C-C_{0}$. \\
\indent We will define three $U$-structures on $C$ as follows.
The core of these structures will be $C_{0}$. For any $k$-ary
basic
relation $R^{{\cal A}_{0}^{+}}$ of ${\cal A}_{0}^{+}$ let \\
\\
\centerline{ $R^{\cal C} = \{s \in {}^{k}C: \ (\exists z \in
R^{{\cal
A}_{0}^{+}}) s \sim z \}$} \\
\\
and for any $k$-ary basic relation $S^{{\cal A}_{1}}$ of ${\cal
A}_{1}$
let \\
\\
\centerline{$S^{{\cal A}_{1}} = \{ s \in {}^{k}C: \ (\exists z \in
S^{{\cal
A}_{1}}) s \sim z \}$.} \\
\\
Finally let \\
\\
\indent ${\cal C}_{0}^{+} = \langle C, C_{0},R^{\cal C} \rangle_{R
\in voc({{\cal A}_{0}^{+}})}$, \indent ${\cal C}_{1} = \langle C,
C_{0},S^{\cal C} \rangle_{S \in voc({\cal A}_{1})}$ \ and \\
\indent let ${\cal C}_{0}$ be the $V$-reduct of ${\cal
C}_{0}^{+}$. \\
\\
By Theorem \ref{equivthm} (1) ${\cal A}_{0}$ and ${\cal A}_{1}$
are ${\cal U}_{n}$-elementary substructures of ${\cal C}_{0}$ and
${\cal C}_{1}$, respectively. Therefore, since ${\cal A}_{0}$ and
${\cal A}_{1}$ are strong $U$-structures, by Theorem
\ref{equivthm} (1) ${\cal C}_{0}$ and ${\cal C}_{1}$ are strong
$U$-structures and moreover ${\cal C}_{0}$ and ${\cal C}_{1}$ are
${\cal U}_{n}$-elementarily equivalent. Similarly, ${\cal
C}_{0}^{+}$ is a model of $\varphi$ (and is a strong
$U$-structure). Let $f$ be the identity function on $C$. By
Theorem \ref{equivthm} (2) $f$ is an isomorphism between ${\cal
C}_{0}$ and ${\cal C}_{1}$. Let ${\cal C}_{1}^{+}$ be the
expansion of ${\cal C}_{1}$ for which $f$ remains an isomorphism
between ${\cal C}_{0}^{+}$ and ${\cal C}_{1}^{+}$ (that is, for
every $R \in voc({\cal C}_{0}^{+})-voc({\cal C}_{0})$ interpret
$R^{{\cal C}_{1}^{+}}$ as $R^{{\cal C}_{1}^{+}} = R^{{\cal
C}_{0}^{+}}$). Clearly, ${\cal C}_{1}^{+} \models \varphi$ and
${\cal C}_{1}^{+}$ is a strong $U$-structure. Let ${\cal
A}_{1}^{+}$ be the substructure of ${\cal C}_{1}^{+}$ generated
by $A_{1}$. Then by Theorem \ref{equivthm} (1) ${\cal A}_{1}^{+}$
is a ${\cal U}_{n}$-elementary substructure of ${\cal C}_{1}^{+}$
and therefore ${\cal A}_{1}^{+} \models \varphi$ (and clearly,
${\cal A}_{1}^{+}$ is a strong $U$-structure by the last sentence
of the statement of Theorem \ref{equivthm} (1)). In addition
${\cal A}_{1}$ is the $V$-reduct of ${\cal A}_{1}^{+}$ and
therefore ${\cal A}_{1} \in K$.
\end{proof}

\begin{theorem}
\label{wcithm} The logic ${\cal U}_{n}$ has the Weak Craig
Interpolation Property.
\end{theorem}

\begin{proof}
Suppose $\varphi$ and $\psi$ are ${\cal U}_{n}$-formulas such
that $\varphi \models \psi$. Let $V$ be the vocabulary consisting
of the relation symbols occurring both in $\varphi$ and in $\psi$.
Let $K_{0}$ be the class of $V$-reducts of models of $\varphi$
and let $K_{1}$ be the class of $V$-reducts of models of $\neg
\psi$. Clearly, $K_{0}$ and $K_{1}$ are closed under
ultraproducts and by Lemma \ref{elemclosedlemma} $K_{0}$ and
$K_{1}$ are closed under ${\cal U}_{n}$-elementary equivalence.
Since $\varphi \models \psi$, it follows that $K_{0}$ and $K_{1}$
are disjoint. Therefore by the Separation Theorem \ref{septhm}
there is a ${\cal U}_{n}$-formula $\vartheta$ (in the common
vocabulary $V$ of $K_{0}$ and $K_{1}$) such that $K_{0} \models
\vartheta$ and $K_{1} \models \neg \vartheta$. But then $\varphi
\models \vartheta$ and $\vartheta \models \psi$, thus $\vartheta$
is the required weak interpolant.
\end{proof}

\begin{theorem}
\label{scithm} If $n \geq 3$ then the logic ${\cal U}_{n}$ doesn't
have the Strong Craig Interpolation Property.
\end{theorem}

\begin{proof}
Let $P$ and $Q$ be two distinct unary relation symbols.
Throughout this proof we will use the vocabulary consisting the
equality symbol, $P$ and $Q$. Let $\varphi(x,y) = P(x)
\Leftrightarrow \neg P(y)$ and let $\psi(x,y,z) = (Q(x)
\Leftrightarrow Q(z)) \vee (Q(y) \Leftrightarrow Q(z))$. \\
\indent First we show that in the class of strong $U$-structures \\
\\
\centerline{ $(*) \indent \models \varphi \Rightarrow \psi$.} \\
\\
To do this assume ${\cal A} \models \varphi[a,b]$ where ${\cal
A}$ is a strong $U$-structure with core $U_{0}$ and $a,b \in A$.
Since $P$ is a unary definable relation of ${\cal A}$, it follows
from Lemma \ref{basiclemma} (1) that $P \in
\{\emptyset,A,U_{0},A-U_{0} \}$. According to our assumption
${\cal A} \models \varphi[a,b]$, either $P = U_{0}$ or
$P=A-U_{0}$. In both cases it follows that exactly one of $\{a,b
\}$ is in $U_{0}$. Similarly, since $Q$ is a unary definable
relation in ${\cal A}$, by Lemma \ref{basiclemma} (1) it follows
that $Q \in \{ \emptyset, A, U_{0}, A-U_{0} \}$. In the first two
cases ${\cal A} \models \psi [a,b,c]$, for any $c \in A$. Now
suppose $Q$ is either $U_{0}$ or $A-U_{0}$. Then exactly one of
$\{ a,b \}$ is in $Q$. Therefore for any $c \in A$ we have ${\cal
A}
\models \psi [a,b,c]$. Thus, $(*)$ is true. \\
\indent Now suppose, seeking a contradiction, that ${\cal U}_{n}$
has the Strong Craig Interpolation Property. Then there exists a
formula $\vartheta$ in which the only relation symbol may be the
equality-symbol such that $\models (\varphi \Rightarrow \vartheta)
\wedge (\vartheta \Rightarrow \psi)$. Now let ${\cal A} = \langle
A,U_{0},P,Q \rangle$ be the strong $U$-structure described in
Theorem \ref{uexistthm} (2). Let $a \in U_{0}, b \in A-U_{0}$,
$a', b' \in U_{0}, a' \not= b', c \in A-U_{0}-\{ b \}$. Then
${\cal A} \models \varphi [a,b,c]$ therefore ${\cal A} \models
\vartheta [a,b,c]$. Observe that there is a permutation $f$ of $A$
with $f(a) = a', f(b) = b', f(c) = c$. Since the only relation
symbol that may occur in $\vartheta$ is the equality, it follows
that ${\cal A} \models \vartheta [f(a),f(b),f(c)]$ and thus
${\cal A} \models \vartheta[a',b',c]$. Therefore, since
$\vartheta$ is a strong interpolant, ${\cal A} \models
\psi[a',b',c]$ would follow, but this contradicts to the choice
of $a',b',c$.
\end{proof}

Let ${\cal U}_{\omega}$ be the logic \\
\indent $\bullet$ whose formulas are that of usual First Order
Logic with $\omega$ many individual variables (but again, the
vocabularies contain relation symbols only) and \\
\indent $\bullet$ whose models are the strong $U$-structures.
\\
Then ${\cal U}_{\omega}$ does not have the Strong Craig
Interpolation Property because the proofs of Lemma
\ref{basiclemma} (1) and
Theorem \ref{scithm} can be repeated in this case, as well. \\
\indent On the other hand ${\cal U}_{\omega}$ still has the Weak
Craig Interpolation Property. To check this, suppose $\varphi$ and
$\psi$ are formulas of ${\cal U}_{\omega}$ such that $\varphi
\models \psi$. Then there exists an $n \in \omega$ for which
$\varphi$ and $\psi$ are formulas of ${\cal U}_{n}$. It is easy
to see that ''$\varphi \models \psi$ in the sense of ${\cal
U}_{\omega}$'' holds if and only if ''$\varphi \models \psi$ in
the sense of ${\cal U}_{n}$''
Hence by Theorem \ref{wcithm} the required interpolant exists in
${\cal U}_{n}$ and consequently in
${\cal U}_{\omega}$ as well. \\
\indent Thus, ${\cal U}_{\omega}$ is an example for a logic with
infinitely many individual variables that has the Weak Craig
Interpolation Property but does not have the Strong Craig
Interpolation Property.

\section{Definability}
\label{defin}

\ \indent The goal of this section is to prove that ${\cal
U}_{n}$ has the Beth Definability Property. For completeness we
start by recalling the relevant definitions.

\begin{definition}
\label{defdef} Let ${\cal L}$ be a logic, let $L \subseteq L^{+}$
be vocabularies for ${\cal L}$ and suppose $R$ is the unique
relation
symbol of $L^{+}$ not occurring in $L$. Suppose $T^{+}$ is a theory in ${\cal L}[L^{+}]$. \\
$\bullet$ We say that $T^{+}$ implicitly defines $R$ over $L$ if
the following holds. If ${\cal A}, {\cal B} \models T^{+}$ and the
$L$-reducts of ${\cal A}$ and ${\cal B}$ are the same (that is,
the identity function on $A$ is an isomorphism between them) then
${\cal A}$ and ${\cal B}$ are the same. \\
$\bullet$ We say that $R$ can be explicitly defined in $T^{+}$
over $L$ if there is a formula of ${\cal L}[L]$ which is
equivalent with $R$
in every model of $T^{+}$. \\
$\bullet$ We say that ${\cal L}$ has the Beth Definability
Property if for any $L,L^{+},T^{+}$ whenever $T^{+}$ implicitly
defines $R$ over $L$ then $R$ can be explicitly defined in
$T^{+}$ over $L$.
\end{definition}

Now we prove a Svenonius-type definability theorem for ${\cal
U}_{n}$. The construction is essentially the same as Theorem
10.5.1 and Corollary 10.5.2 of \cite{Hodges}.

\begin{theorem}
\label{sventhm} Suppose $L \subseteq L^{+}$ are vocabularies for
${\cal U}_{n}$, $R$ is the unique relation symbol of $L^{+}$ not
occurring in $L$ and $T^{+}$ is a complete theory in $L^{+}$.
Then the
following are equivalent. \\
\indent (1) $R$ can be explicitly defined in $T^{+}$ over $L$. \\
\indent (2) If ${\cal A} \models T^{+}$ and $f$ is an
automorphism of ${\cal A}|_{L}$ then $f$ preserves $R^{\cal A}$ as
well.
\end{theorem}

\begin{proof}
Clearly, (2) follows from (1). To prove the converse implication
suppose (2) holds and suppose, seeking a contradiction, that $R$
cannot be explicitly defined in $T^{+}$ over $L$. Suppose that $R$
is $k$-ary for some $k \leq n$. Expand $L^{+}$ by two $k$-tuples
$\bar{c}, \bar{d}$ which are new constant symbols and let $\Gamma
= \{ \varphi(\bar{c}) \Leftrightarrow \varphi(\bar{d}): \varphi$
is a ${\cal U}_{n}[L]$-formula $\}$. Consider the following first
order theory $\Sigma$ (since in ${\cal U}_{n}$ constant symbols
are not part of the language,
strictly speaking the following $\Sigma$ is not a theory in ${\cal U}_{n}$). \\
\\
\centerline{ $\Sigma = T^{+} \cup \Gamma
\cup \{ R(\bar{c}), \neg R(\bar{d}) \}$. } \\
\\
We claim that every finite subset $\Sigma_{0}$ of $\Sigma$ has a
model whose $L^{+}$-reduct is a strong $U$-structure. To show
this suppose, seeking a contradiction, that $\Sigma_{0}$ is a
finite subset of $\Sigma$ which doesn't have such a model. Let
$\Gamma' = \Sigma_{0} \cap \Gamma = \{ \varphi_{i}(\bar{c})
\Leftrightarrow
\varphi_{i}(\bar{d}): i < m \}$. Observe that \\
\\
\centerline{$(*)$ if ${\cal A} \models T^{+}$ is
a strong $U$-structure, $\bar{a},\bar{b} \in A$, $\langle {\cal
A}, \bar{a},\bar{b} \rangle \models \Gamma'$ and ${\cal A} \models
R(\bar{b})$} \\
\centerline{then ${\cal A} \models R(\bar{a})$ } \\
\\
because otherwise $\langle {\cal A}, \bar{b},\bar{a} \rangle$
would be a model of $\Sigma_{0}$ whose $L^{+}$-reduct is a strong
$U$-structure. Let $\Phi = \{ \varphi_{i}(\bar{v}): i < m \}$.
Suppose ${\cal A} \models T^{+}$ and $\bar{a} \in A$. Then the
$\Phi$-type of $\bar{a}$ in ${\cal A}$ is defined as follows: \\
\\
\centerline{ $\Phi-tp^{\cal A}(\bar{a}) = \{ \varphi_{i}(\bar{v}):
{\cal A} \models \varphi_{i}(\bar{a}), i < m \} \cup \{ \neg
\varphi_{j}(\bar{v}): {\cal A} \not\models \varphi_{j}(\bar{a}),
i < m \}$.} \\
\\
Let $\varrho = \bigvee \{ \bigwedge \psi: $ there are a strong
$U$-structure ${\cal A} \models T^{+}$ and $\bar{a} \in R^{\cal
A}$ such that $\psi = \Phi-tp^{\cal A}(\bar{a}) \}$. Clearly,
$\varrho$ is a formula of ${\cal U}_{n}[L]$. We claim that
$\varrho$ defines explicitly $R$ in $T^{+}$ over $L$. To verify
this suppose ${\cal A} \models T^{+}$. If $\bar{a} \in R^{\cal
A}$ then $\bigwedge (\Phi-tp^{\cal A}(\bar{a}))$ is a disjunctive
component of $\varrho$ therefore ${\cal A} \models
\varrho(\bar{a})$. Thus, the relation defined by $\varrho$ in
${\cal A}$ contains $R^{\cal A}$. Conversely, suppose ${\cal A}
\models \varrho(\bar{b})$. Then there is a disjunctive component
$\bigwedge \psi$ of $\varrho$ such that ${\cal A} \models
\bigwedge \psi(\bar{b})$ and there are another strong
$U$-structure ${\cal A}' \models T^{+}$ and $\bar{a}' \in
R^{{\cal A}'}$ such that $\bigwedge \psi = \bigwedge (
\Phi-tp^{{\cal A}'}(\bar{a}'))$. Thus, ${\cal A}' \models \exists
\bar{v} (R(\bar{v}) \wedge \bigwedge \psi(\bar{v}))$. This last
formula is a ${\cal U}_{n}$-formula, and since $T^{+}$ is
complete, ${\cal A} \models \exists \bar{v} (R(\bar{v}) \wedge
\bigwedge \psi(\bar{v}))$. Thus, there is $\bar{a} \in R^{\cal A}$
such that $\bigwedge(\Phi-tp^{\cal A} (\bar{a})) = \bigwedge \psi
= \bigwedge (\Phi-tp^{\cal A}(\bar{b}))$. Therefore by
$(*)$ it follows that ${\cal A} \models R(\bar{b})$. \\
\indent We proved that $\varrho$ explicitly defines $R$ in
$T^{+}$ over $L$. This is impossible because we assumed that $R$
cannot be explicitly defined. Hence every finite subset of
$\Sigma$ has a model whose $L^{+}$-reduct is a strong $U$-structure.\\
\indent Let $\langle {\cal A}, \bar{a}, \bar{b} \rangle$ be an
ultraproduct of the above models of finite subsets of $\Sigma$
for which $\langle {\cal A}, \bar{a}, \bar{b} \rangle \models
\Sigma$. By Lemma \ref{upclosed} the $L^{+}$-reduct of it (which
is ${\cal A}$) is a strong $U$-structure. Since $\Gamma \subseteq
\Sigma$, it follows that $\bar{a}$ and $\bar{b}$ satisfies the
same ${\cal U}_{n}[L]$-formulas. Therefore by Theorem
\ref{autothm} there is an automorphism of the $L$-reduct of
${\cal A}$ mapping $\bar{a}$ onto $\bar{b}$. This automorphism
doesn't preserve $R^{\cal A}$, contradicting to (2). This proves
that $R$ can be explicitly defined in $T^{+}$ over $L$.
\end{proof}

\begin{theorem}
\label{bdthm} The logic ${\cal U}_{n}$ has the Beth Definability
Property.
\end{theorem}

\begin{proof}
Let $L,L^{+}, R$ and $T^{+}$ be as in Definition \ref{defdef} and
assume $T^{+}$ implicitly defines $R$ over $L$. We have to show
that $R$ can be explicitly defined in $T^{+}$ over $L$. \\
\indent First suppose that $T^{+}$ is a complete theory.
%
%
Suppose ${\cal A}$ is a model of $T^{+}$ and $f$ is an
automorphism of ${\cal A}|_{L}$. We claim that $f$ preserves
$R^{\cal A}$ as well. To see this, define another structure
${\cal B}$ as follows. The universe of ${\cal B}$ is $A$. For any
subset $X$ of (a direct power of) $A$ the $f$-image of $X$ will
be denoted by $f[X]$. For every $P \in L$ let $P^{\cal B} =
f[P^{\cal A}]$, let $R^{\cal B} = f[R^{\cal A}]$ and let
$U'=f[U]$ where $U$ is the core of ${\cal A}$. Since $f$ is an
automorphism of ${\cal A}|_{L}$, it follows that ${\cal A}|_{L} =
{\cal B}|_{L}$. In addition, $f$ is an isomorphism between
$\langle {\cal A},U \rangle$ and $\langle {\cal B},U' \rangle$.
Therefore ${\cal B}$ is a strong $U$-structure with core $U'$ and
${\cal B} \models T^{+}$.
%
Since $T^{+}$ implicitly defines $R$ over $L$, it follows that
$R^{\cal A} = R^{\cal B}$, that is, $f$ preserves $R^{\cal A}$.
Since ${\cal A}$ and $f$ were chosen arbitrarily, it follows that
every automorphism of the $L$-reduct of a model of $T^{+}$ also
preserves the interpretation of $R$. Therefore by Theorem
\ref{sventhm} $R$
can be explicitly defined in $T^{+}$ over $L$. \\
\indent Now let $T^{+}$ be an arbitrary (not necessarily
complete) theory which implicitly defines $R$ over $L$. We claim
that there is a finite set $\Phi = \{
\varphi_{0},...,\varphi_{m-1} \}$ of ${\cal U}_{n}[L]$-formulas
such that if ${\cal A} \models
T^{+}$ then \\
\\
\centerline{$(*) \indent {\cal A} \models \bigvee_{i < m} (\forall
v_{0}...\forall v_{n-1} (R \Leftrightarrow \varphi_{i} ))$.} \\
\\
For if not, then for any finite set $\Phi$ of ${\cal
U}_{n}[L]$-formulas it would exist a model of $T^{+}$ in which $R$
would be different from all the relations defined by the members
of $\Phi$. Forming an ultraproduct of these models it would exist
a strong $U$-structure ${\cal A} \models T^{+}$ in which $R^{\cal
A}$ would not be definable in ${\cal A}|_{L}$. But then $T' = \{
\varphi: {\cal A} \models \varphi, \varphi $ is a ${\cal
U}_{n}[L^{+}]$-formula $\}$ would be a complete theory and since
$T^{+} \subseteq T'$, $T'$ also implicitly defines $R$ over $L$.
Therefore by the second paragraph of this proof $R$ would be
explicitly definable in $T'$ and particularly, there would be a
${\cal U}_{n}[L]$-formula which would define $R^{\cal A}$ in
${\cal A}$; a
contradiction. Therefore $(*)$ is established. \\
\indent Now for each ${\cal A} \models T^{+}$ let $\nu({\cal A})$
be the smallest $i \in m$ for which ${\cal A} \models R
\Leftrightarrow \varphi_{i}$ and let $K_{i} = \{ {\cal A}|_{L}:
{\cal A} \models T^{+}, \nu({\cal A}) = i \}$. Clearly, the
classes $K_{i}$ are pairwise disjoint and closed under
ultraproducts. In fact they are closed under ${\cal
U}_{n}$-elementary equivalence because of the following. Suppose
${\cal A} \in K_{i}$ and ${\cal A}$ and ${\cal B}$ are ${\cal
U}_{n}$-elementarily equivalent. Let $U_{0}, V_{0}$ be the cores
of ${\cal A}$ and ${\cal B}$, respectively. Let $C_{0} \subseteq
C$ be two sets such that $|C_{0}| \geq |U_{0}|,|V_{0}|$ and
$|C-C_{0}| \geq |A-U_{0}|,|B-V_{0}|$. Then we may assume that
$U_{0}, V_{0} \subseteq C_{0}, A-U_{0},B-V_{0} \subseteq
C-C_{0}$. We will define two $U$-structures on $C$ as follows. If
$R$ is any $m$-ary relation symbol in $L$ then let $R^{{\cal
C}_{0}} = \{ s \in {}^{m}C: \ (\exists z \in R^{\cal A}) s \sim z
\}$ and let $R^{{\cal C}_{1}} = \{ s \in {}^{m}C: \ (\exists z \in
R^{\cal B}) s \sim z \}$. Then by Theorem \ref{equivthm} (1)
${\cal A}$ and ${\cal B}$ are ${\cal U}_{n}$-elementary
substructures of ${\cal C}_{0}$ and ${\cal C}_{1}$, respectively.
Therefore ${\cal C}_{0}$ and ${\cal C}_{1}$ are strong
$U$-structures and ${\cal U}_{n}$-elementarily equivalent with
each other. Hence by Theorem \ref{equivthm} (2) the identity
function on $C$ is an isomorphism between ${\cal C}_{0}$  and
${\cal C}_{1}$. Since ${\cal A} \in K_{i}$, $i$ is the smallest
number for which $R$ and $\varphi_{i}$ are equivalent in ${\cal
A}$. Let $R^{{\cal C}_{0}}$ be the relation defined by
$\varphi_{i}$ in ${\cal C}_{0}$. Since ${\cal A}$ is a ${\cal
U}_{n}$-elementary
substructure of ${\cal C}_{0}$, it follows that \\
\\
\indent (i) \indent $\langle {\cal C}_{0}, R^{{\cal C}_{0}}
\rangle \models T^{+}$. \\
\\
Since $T^{+}$ implicitly defines $R$, this is the only way to
extend ${\cal C}_{0}$ to a model of $T^{+}$. In particular, \\
\\
\indent (ii) \indent for every $j < i$ we have ${\cal C}_{0} \not
\models R^{{\cal C}_{0}} \Leftrightarrow \varphi_{j}$. \\
\\
Since the identity function of $C$ is an isomorphism between
${\cal C}_{0}$ and ${\cal C}_{1}$, the above (i) and (ii) are
true for ${\cal C}_{1}$ as well. Let $R^{\cal B} = R^{{\cal
C}_{0}} \cap {}^{m} B $. Then by Theorem \ref{equivthm} (1)
$\langle {\cal B}, R^{\cal B} \rangle \models T^{+}$ and $i$ is a
smallest number for which $\varphi_{i}$ defines $R^{\cal B}$ in
${\cal
B}$. Thus, ${\cal B} \in K_{i}$, as desired. \\
\indent Now by Theorem \ref{septhm} for every $i \in m$ there is a
${\cal U}_{n}[L]$-formula $\varrho_{i}$ such that $K_{i} \models
\varrho_{i}$ and $\cup_{j \in m - \{ i \}} K_{i} \models \neg
\varrho_{i}$. Finally let $\psi = \bigvee_{i \in m} (\varrho_{i}
\wedge \varphi_{i})$. It is easy to check that $\psi$ is
equivalent with $R$ in every model of $T^{+}$, thus $R$ can be
explicitly defined in $T^{+}$ over $L$, as desired.
\end{proof}

\section{Cylindric Algebraic Consequences}
\label{cyla}

\ \indent By translating the results of the previous sections to
Algebraic Logic, in this section we prove that for finite $n \geq
3$, there is a (finitely axiomatizable) subvariety of
$RCA_{n}$ that has the Strong Amalgamation Property but doesn't
have the Superamalgamation Property (the definitions of these
properties can be found for example in \cite{ans} before
Definition 6.14). As we mentioned this settles a problem
of Pigozzi in \cite{pigo}. \\
\indent We assume that the reader is familiar with the theory of
cylindric algebras. Some basic facts on this topic have been
recalled before Lemma \ref{upclosed}. For more details we refer to
\cite{HMT1} and \cite{HMT2}. \\
%
%
\indent If $K$ is a class of algebras then ${\bf S} K$ and ${\bf
P} K$ denote the classes of (isomorphic copies of) subalgebras of
members of $K$ and (isomorphic copies of) direct products of
members of $K$, respectively. Similarly, ${\bf Up} K$ denotes the
class of (isomorphic copies of) ultraproducts of members of $K$.
For other algebraic notions and notation we refer to \cite{bs}.

\begin{definition} $US_{n}$ and $U_{n}$ are defined to be the following
subclasses of $RCA_{n}$: \\
\\
\indent $US_{n} =  \{ Cs_{n}({\cal A}): \ {\cal A}$ is a strong
$U$-structure $ \}$. \\
\indent $U_{n} = {\bf SP} US_{n}$.
\end{definition}

\begin{theorem}
\label{varthm}
$U_{n}$ is a finitely axiomatizable variety.
\end{theorem}

\begin{proof}
Lemma \ref{basiclemma} (3) implies that there is a natural number
$N_{1}$ such that for all strong $U$-structure ${\cal A}$ we have
$|Cs_{n}({\cal A})| \leq N_{1}$ (we already observed this in the
proof of Lemma \ref{upclosed}). Therefore $US_{n}$ is finite and
hence ${\bf Up}US_{n} = US_{n}$.
So $U_{n} = {\bf SP}US_{n} \subseteq {\bf SPUp}US_{n} = {\bf SP}
US_{n} = U_{n}$. Hence $U_{n}$ is the quasi-variety generated by
$US_{n}$. The cylindric term $c_{0}...c_{n-1}(x)$ is a
switching-function in ${\bf S}US_{n}$ therefore the quasi-variety
and the variety generated by
$US_{n}$ coincide. Thus $U_{n}$ is the variety generated by $US_{n}$. \\
%
%
%
\indent Finally observe that $U_{n}$ is congruence-distributive
since it has a Boolean reduct. Thus, $U_{n}$ is a finitely
generated congruence-distributive variety and hence by Baker's
Theorem it is finitely axiomatizable (see \cite{baker} or
\cite{bs}).
\end{proof}

Now we return to the choice of naming our logic ${\cal U}_{n}$ and
the classes $US_{n}$ and $U_{n}$. By Lemma \ref{basiclemma} (2)
every member of $US_{n}$ is a subalgebra of the $n$-dimensional
Cylindric Set Algebra generated by one UNARY relation: by the core
of the corresponding structure. So
''$U$'' stands for ''unary''. \\
\indent
Now we are ready to prove the main theorem of the paper.

\begin{theorem}
\label{mainthm}
(1) $U_{n}$ has the Strong Amalgamation Property. \\
\indent (2) $U_{n}$ doesn't have the Superamalgamation Property,
if $n \in \omega, n \geq 3$.
\end{theorem}

\begin{proof}
(1) By theorem \ref{bdthm} ${\cal U}_{n}$ has the Beth
Definability Property and therefore by \cite{nemeti} the
epimorphisms of $U_{n}$ are surjective (see also \cite{ans},
Theorem 6.11). By Theorem \ref{wcithm} ${\cal U}_{n}$ has the Weak
Craig Interpolation Property and by Theorem \ref{varthm} $U_{n}$
is a variety. Thus, by Theorem 6.15(i) of \cite{ans} (see also the
beginning of Section 7 therein) $U_{n}$ has the Amalgamation
Property. Since $U_{n}$ is a variety, it follows from
Propositions 1.9 and 1.11 of \cite{KMPT} (see also Proposition 6.3
therein)
that $U_{n}$ indeed has the Strong Amalgamation Property. \\
\indent (2) By Theorem \ref{scithm} ${\cal U}_{n}$ doesn't have
the Strong Craig Interpolation Property and therefore by Theorem
6.15 (ii) of \cite{ans} $U_{n}$ doesn't have the
Superamalgamation Property.
\end{proof}

{\bf Acknowledgement.} Thanks are due to Alice Leonhardt for the
beautiful typing of an early and preliminary version of this
notes.

\bigbreak
\leftline{Alfr\'ed R\'enyi Institute of Mathematics}
\leftline{Hungarian Academy of Sciences}
\leftline{Budapest Pf. 127}
\leftline{H-1364 Hungary}
\leftline{sagi@renyi.hu}

\ \\

\bigbreak
\leftline{Department of Mathematics}
\leftline{Hebrew University}
\leftline{91904 Jerusalem, Israel}
\leftline{shelah@math.huji.ac.il}


\begin{thebibliography}{99}

\bibitem {ans} {\sc H. Andr\'eka, I. N\'emeti, I. Sain,}
{\it Algebraic Logic, \/} in Handbook of Philosophical Logic
(eds. D. M. Gabbay and F. Guenthner), 2nd edition, Kluwer
Academic Publishers, (2001).
\bibitem {baker}{\sc K. Baker,}
{\it Finite Equational Bases for Finite Algebras in a
Congruence-Distrubutive Equational Class, \/}, Advances in
Mathematics 24, pp.204-243, (1977).
\bibitem {bs} {\sc S. Burris, H. P. Sankappanavar,}
{\it A Course in Universal Algebra, \/} Spinger Verlag, New York
(1981).
\bibitem {chk} {\sc C.C. Chang, H.J. Keisler,}
{\it Model Theory, \/} North--Holland, Amsterdam (1973).
\bibitem {comer} {\sc Comer,}
{\it Classes without the Amalgamation Property, \/} Pacific J.
Math. 28, pp. 309-318, (1969).
\bibitem {HMT1} {\sc L. Henkin, J. D. Monk, A. Tarski,}
{\it Cylindric Algebras Part 1, \/} North-Holland, Amsterdam
(1971).
\bibitem {HMT2} {\sc L. Henkin, J. D. Monk, A. Tarski,}
{\it Cylindric Algebras Part 2, \/} North-Holland, Amsterdam
(1985).
\bibitem {Hodges} {\sc W. Hodges,}
{\it Model theory, \/} Cambridge University Press, (1997).
\bibitem {KMPT} {\sc E. W. Kiss, L. M\'arki, P. Pr\H{o}hle and W. Tholen,}
{\it Categorical Algebraic Properties. A Compendium on
Amalgamation, Congruence Extension, Epimorphisms, Residual
Smallness and Injectivity \/} Studia Sci. Math. Hungarica 18, pp.
79-141, (1983).
\bibitem {maksi} {\sc L. Maksimova,}
{\it Beth's Property, Interpolation and Amalgamation in Varieties
of Modal Algebras, \/} (Russian) Doklady Akademii Nauk. SSSR. vol.
319 (1991)  no. 6, pp. 1309-1312.
\bibitem {nemeti} {\sc I. N\'emeti,}
{\it Beth Definability Property is Equivalent with Surjectiveness
of Epis in general Algebraic Logic, \/} Tehchnical Report of the
Mathematical Institute of Hungarian Academy of Sciences,
Budapest, 1983.
\bibitem {pigo} {\sc D. Pigozzi,}
{\it Amalgamation, Congruence Extension and Interpolation
Properties in Algebras, \/} Algebra Universalis Vol. 1 No. 3.,
pp. 269-349, (1972).
\bibitem {clasth} {\sc S. Shelah,}
{\it Classification theory, \/} North--Holland, Amsterdam (1990).
\end{thebibliography}
\end{document}